\documentclass[a4paper, 12pt]{article}

\usepackage{amsmath, amstext, amsxtra,amsthm, amscd, amsgen, amsbsy, amsopn, amsfonts, latexsym,
amssymb,amscd}

\begin{document}


\def\seminegative{semi-negative}
\def\LFF{Lefschetz fixed point formula}

\def\tind{t\text{-}ind}

\def\proofend{\hbox to 1em{\hss}\hfill $\blacksquare $\bigskip }

\newtheorem{theorem}{Theorem}[section]
\newtheorem{proposition}[theorem]{Proposition}
\newtheorem{lemma}[theorem]{Lemma}
\newtheorem{remarks}[theorem]{Remarks}
\newtheorem{remark}[theorem]{Remark}
\newtheorem{definition}[theorem]{Definition}
\newtheorem{corollary}[theorem]{Corollary}
\newtheorem{example}[theorem]{Example}
\newtheorem{assumption}[theorem]{Assumption}
\newtheorem{problem}[theorem]{Problem}
\newtheorem{question}[theorem]{Question}
\newtheorem{conjecture}[theorem]{Conjecture}

\def\fib{fibre}
\def\vsp{vector space}
\def\b{bun\-dle}
\def\pb{principal \b }
\def\vb{vector \b }
\def\fb{\fib\ \b }
\def\tb{tangent\ \b }
\def\nb{normal\ \b }

\def\mfd{manifold}

\def\Z{{\mathbb Z}}
\def\R{{\mathbb R}}
\def\Q{{\mathbb Q}}
\def\C{{\mathbb C}}
\def\N{{\mathbb N}}
\def\H{{\mathbb H}}
\def\Zp #1{{\mathbb Z }/#1{\mathbb Z}}
\def\cpt{compact}

\def\paperref#1#2#3#4#5#6{\text{#1:} #2, {\em #3} {\bf#4} (#5)#6}
\def\bookref#1#2#3#4#5#6{\text{#1:} {\em #2}, #3 #4 #5#6}


\title{Bordism-finiteness and semi-simple group actions}
\author{Anand Dessai}
\date{}
\maketitle

\begin{abstract}
\noindent
We give bordism-finiteness results for smooth $S^3$-manifolds. Consider the class of oriented manifolds which admit an 
$S^1$-action with isolated fixed points such that the action extends 
to an $S^3$-action with fixed point. We exhibit various subclasses, 
characterized by an upper bound for the Euler characteristic and 
properties of the first Pontrjagin class $p_1$, for example $p_1=0$, 
which contain only finitely many oriented bordism types in any given 
dimension. Also we show finiteness results for homotopy 
complex projective spaces and complete intersections with $S^3$-action 
as above. 
\end{abstract}

\section{Introduction}\label{intro} A direct consequence of the Atiyah-Bott-Segal-Singer fixed point theorem
in equivariant index theory is the fact that all the Pontrjagin numbers of an oriented closed manifold 
with smooth fixed point free $S^1$-action vanish. Since the oriented bordism ring is determined by Pontrjagin and 
Stiefel-Whitney numbers this result may be interpreted as a 
bordism-finiteness theorem: In a given dimension the class of oriented 
closed manifolds which support a fixed point free $S^1$-action 
contains only finitely many oriented bordism types. 

In this paper we give bordism-finiteness results for smooth actions on 
connected closed oriented manifolds with possible non-empty fixed point set. We 
also look at the related question of diffeomorphism-finiteness in a 
given homotopy type. 

In Riemannian geometry finiteness theorems (involving bounds for 
curvature, diameter, volume 
 etc.) go back to the 
work of Cheeger (cf. \cite{Ch}) and are the topic of active research 
since then. It is natural to look what happens if the geometrical 
bounds arise from symmetries, i.e. group actions. 

Concerning the questions on bordism- and diffeomorphism-finiteness it 
is desirable to obtain positive results for groups as ``small'' as 
possible. Our results involve smooth actions by the semi-simple group 
$S^3$. Before we state these we take a brief look at actions by 
``smaller'' groups. 

For finite groups one cannot expect any general finiteness results, 
even if one restricts to free actions or to a fixed homotopy type. 
Dovermann and Masuda used equivariant surgery to construct infinitely 
many homotopy $\C P^3$'s with effective smooth $\Z _p$-action for 
every prime number $p$ (cf. \cite{DoMa} and the references therein). 
L\"offler and Rau\ss en (cf. \cite{LoRa}) used Sullivan's theory of 
minimal models to construct free $\Z 
_p$-actions on certain high connected manifolds for any large enough 
prime number $p$ and conjectured that non-trivial $\Z _p$-actions 
exist on any simply-connected manifold (cf. also \cite{We}; for the 
question of bordism-finiteness cf. \cite{CoFl} and \cite{toDi}). On 
the other hand, by a result of Schultz, one knows that in dimension 
$\geq 4$ any oriented bordism class contains infinitely many manifolds 
with infinite fundamental group which admit no effective finite group 
action (cf. \cite{Sc}). 

Next we discuss the question of bordism-finiteness for circle actions. 
By taking products or equivariant connected sum along invariant 
submanifolds one can construct new $S^1$-manifolds from given pieces. 
Using these methods it is easy to show that for any given oriented 
manifold a suitable multiple (disjoint union) is bordant to a 
connected oriented manifold with non-trivial $S^1$-action. On the 
other hand the 
\LFF \ in equivariant index theory (cf. \cite{AtSe}, \cite{AtSiIII}, \cite{AtBo}, \cite{Bo}) implies that 
the bordism class of a manifold with $S^1$-action is modulo torsion 
completely determined by the local geometry of the action near the 
fixed point set (see Cor. \ref{localgeometry}). 

In order to get interesting bordism-finiteness results one needs 
additional conditions which limit the constructions mentioned above 
but are not too restrictive on the local geometry of the action. This 
leads us to consider $S^1$-actions with {\it isolated} fixed points 
which satisfy a prescribed upper bound for the number of fixed points 
(by the classical \LFF \ the number of $S^1$-fixed points is just the 
Euler characteristic). In this way we exclude products with trivial 
$S^1$-manifolds and limit the use of the connected sum construction 
but still allow various local geometries although the non-equivariant 
part of these is trivial.

For these $S^1$-actions bordism-finiteness still fails (see the 
examples in Section \ref{examples}). Taking products one concludes the 
same for torus actions as long as the dimension of the torus is small 
compared to the dimension of the manifold. We note that the absence of 
strong implications on the bordism type for torus actions is 
underlined by the following result of Buchstaber and Ray (cf. 
\cite{BuRa}): The complex bordism ring is generated by toric 
manifolds. In particular, any $2n$-dimensional stable almost complex 
manifold is bordant to a manifold which admits an action by the 
$n$-dimensional torus. 

Since bordism-finiteness fails for abelian actions it is natural to 
look at semi-simple group actions next. Our results involve actions of 
$S^3$ with the following property: 
\begin{equation}S^3\text{ acts with fixed point and the }S^1\text{-action has isolated
fixed points}\tag{$\ast $}\end{equation} for some (and hence every) 
fixed subgroup $S^1\hookrightarrow S^3$. In turns out that for such 
actions the question whether bordism-finiteness holds depends on 
properties of the first Pontrjagin class. A special case is the 
following 

\begin{theorem}\label{p1zero} Let $C$ and $m$ be natural numbers. The class 
of connected $m$-dimensional oriented manifolds with vanishing first 
Pontrjagin class and Euler characteristic $\leq C$ which admit an 
$S^3$-action satisfying $(\ast )$ contains only finitely many oriented 
bordism types. 
\end{theorem} 

\noindent The conclusion holds in the more 
general situation where a negative multiple of the first Pontrjagin 
class $p_1$ is a sum of squares (see Th. \ref{theobordfin}). Moreover 
bordism-finiteness holds if one refines the bordism ring taking into 
account the condition on $p_1$ (see Th. \ref{bkfinite}). To prove 
Theorem \ref{p1zero} and its refinements we combine the \LFF \ in 
equivariant index theory with the following well known properties of 
$S^3$ (which hold for any semi-simple group): In any fixed dimension 
the group $S^3$ admits only finitely many non-equivalent 
representations. Any one-dimensional $S^3$-representation is trivial. 

These properties imply that the local geometry of the induced 
$S^1$-action at the $S^3$-fixed point is determined up to finite 
ambiguity. The condition on the first Pontrjagin class guarantees that 
this is also true at the other $S^1$-fixed points. Finally one applies 
the 
\LFF \ in equivariant index theory to conclude that the oriented bordism type (resp. its 
refinement) is determined up to finite ambiguity. 

 Regarding the problem of 
diffeomorphism-finiteness in a given homotopy type we use similar 
arguments to show 
\begin{theorem}\label{hcp} In a fixed dimension $2n$ there are only 
finitely many homotopy complex projective spaces with an $S^3$-action 
satisfying 
 $(\ast )$.\end{theorem}

\noindent
By simply-connected surgery theory one knows that in any dimension 
$2n\geq 6$ there are infinitely many homotopy complex projective 
spaces distinguished by their Pontrjagin classes (cf. \cite{Hs}). 
According to Theorem \ref{hcp} almost all of them do not admit an 
$S^3$-action satisfying 
 $(\ast )$. It is interesting to compare the result above with a conjecture
 of Petrie (cf. \cite{Pe}) which states that the total Pontrjagin class of a homotopy $\C P^n$ with 
non-trivial $S^1$-action is standard. The conjecture is known to be 
true if $n\leq 4$ (cf. \cite{Ja}; for related results cf. 
for example the survey \cite{Do} as well as \cite{DAN} and the 
references therein). It implies diffeomorphism-finiteness for homotopy 
complex projective spaces with non-trivial $S^1$-action. 

As indicated before our methods are rather classical, involving index 
theory and properties of semi-simple groups, and work best if a 
negative multiple of the first Pontrjagin class is a sum of squares. 
On the other hand the theory of elliptic genera led to applications 
for $Spin^c$-manifolds with $S^3$-action if the first Pontrjagin class 
is a sum of squares (cf. \cite{DAN}, \cite{DHCP}, \cite{Li}). 
Combining both methods one obtains further information on manifolds 
with $S^3$-action. 

This paper is structured in the following way. Section \ref{weights} 
deals with tangential weights and weights of equivariant complex line 
bundles for manifolds with admit an $S^3$-action with fixed point. We 
give conditions in terms of the first Pontrjagin class which guarantee 
that these weights are determined up to finite ambiguity. In Section 
\ref{bordfin} we extend Theorem \ref{p1zero} to the case that a 
negative multiple of the first Pontrjagin class is a sum of squares 
(see Theorem \ref{theobordfin}). To prove this theorem and related 
results we use the 
\LFF 
\ in equivariant index theory. In Section \ref{examples} we show, by example, that 
the bordism-finiteness result given in Theorem \ref{theobordfin} fails 
if one relaxes the conditions on the $S^3$-action or the condition on 
the first Pontrjagin class. In Section \ref{homcp} we prove the result 
on homotopy complex projective spaces (see Theorem \ref{hcp} above). 
We use similar methods to show that in a given dimension there are 
only finitely many complete intersections which admit an $S^3$-action 
satisfying $(\ast )$. 

\section{Weights}\label{weights}
In this section we give some information on the tangential weights and 
the weights of complex line bundles over a $2n$-dimensional manifold 
$M$ which supports a non-trivial $S^3$-action with fixed point. 

By a result of Hattori and Yoshida (cf. \cite{HaYo}) the $S^3$-action 
lifts uniquely to any complex line bundle $L\to M$ and we shall always 
do so. Fix $S^1\hookrightarrow S^3$. At an $S^1$-fixed point $L$ 
reduces to a complex one-dimensional $S^1$-representa\-tion with 
character $\lambda 
\mapsto 
\lambda ^{a}$. We call $a$ the {\bf weight} of $L$ at this point. Note 
that at the $S^3$-fixed point the weight has to vanish since any 
complex one-dimensional $S^3$-representation is trivial (one way to 
characterize semi-simple compact Lie groups). 

At an $S^1$-fixed point the tangent bundle $TM$ reduces to a real 
$S^1$-represen\-ta\-tion which we identify (non-canonically) with a 
complex representation with character $\lambda 
\mapsto \sum 
\lambda^{m_i}$. We call $m_1,\ldots ,m_n$ the {\bf tangential weights} of $M$ at 
the fixed point. Note that the tangential weights are only well 
defined up to sign. At an $S^3$-fixed point $pt$ the 
$S^1$-representation $TM_{\vert pt}$ is induced from an 
$S^3$-representation. Since there are only finitely many equivalence 
classes of such representations in a given dimension (another way to 
characterize semi-simple compact Lie groups) the tangential weights at 
$pt$ are determined by the dimension of $M$ up to finite ambiguity, 
i.e. the set of tangential weights at $pt$ belongs to a finite set 
which only depends on $2n$. The next lemma is used in the following 
section to derive the bordism-finiteness results mentioned in the 
introduction. 
\begin{lemma}\label{mainlemma} Let $L_1,\ldots ,L_k$ be $S^3$-equivariant 
complex line bundles over $M$ and $N$ a positive integer. Assume the 
first Pontrjagin class of the bundle $E:=(L_1+\ldots +L_k )+N\cdot TM$ 
vanishes. Then the tangential weights at any $S^1$-fixed point are 
determined by the dimension $2n$ of $M$ up to finite ambiguity. The 
weights of $L_i$ are determined by $N$ and $2n$ up to finite 
ambiguity. 
\end{lemma}

\bigskip
\noindent
{\bf Proof:} Note that the weights of $E$ at an $S^1$-fixed point are 
just the weights of the line bundles $L_j$ and the tangential weights 
taken with multiplicity $N$. Let $Y_1,\ldots ,Y_{C^\prime }$ denote 
the connected components of the fixed point manifold $M^{S^1}$. We 
assume that the $S^3$-fixed point $pt$ belongs to $Y_1$. Let 
$m_{s,1},\ldots ,m_{s,n}$ denote the tangential weights and let 
$a_{s,j}$ denote the weight of $L_j$ at $Y_s$. 

Since the $S^3$-equivariant vector bundle $E$ has vanishing first 
Pontrjagin class it follows from a spectral sequence argument (see the 
lemma below) that the weights of $E$ at $Y_s$ satisfy 
\begin{equation}\label{boundweight}
\sum_{j=1}^k a_{s,j}^2+N\cdot \sum _{i=1}^n 
m_{s,i}^2=C^{\prime 
\prime },\end{equation}
 where $C^{\prime \prime }$ is a
constant which does not depend on $Y_s$. At the $S^3$-fixed point $pt 
\in Y_1$ the weights of $L_j$ vanish and the tangential weights 
$m_{1,i}$ are determined by the dimension of $M$ up to finite 
ambiguity. Hence $C^{\prime \prime }=N\cdot \sum 
_{i=1}^n m_{1,i}^2$ is bounded from above by $N$ times a constant which only depends 
on $2n$. By formula (\ref{boundweight}) the tangential weights (resp. 
the weights of $L_j$) at any $S^1$-fixed point are determined by $2n$ 
(resp. $2n$ and $N$) up to finite ambiguity.\proofend 

\begin{lemma}\label{specseq}
Let $E\to M$ be a $2r$-dimensional $S^3$-equivariant vector bundle 
with weights $e_{s,1},\ldots ,e_{s,r}$ at the $S^1$-fixed point 
component $Y_s$. If $p_1(E)=0$ then $\sum _{i=1}^r e_{s,i}^2$ is 
independent of $Y_s$. 
\end{lemma}

\bigskip
\noindent
{\bf Proof:} We use the Leray-Serre spectral sequence for the Borel 
construction of $E\to M$ (for details cf. \cite{DDISS}, Lemma A.6). 
Let $M_{S^3}:=ES^3\times _{S^3}M$ and let $\pi :M_{S^3}\to BS^3$ be 
the fibration with fibre $i:M\hookrightarrow M_{S^3}$ associated to 
the universal $S^3$-principal bundle $ES^3\to BS^3$ and the 
$S^3$-space $M$. Since $H^*(BS^3;\Z )$ is concentrated in degree $4*$ 
it follows from the Leray-Serre spectral sequence that the sequence 
\begin{equation}\label{eqspecseq} H^4(BS^3;\Z )\overset
 {\pi ^*}\to H^4(M_{S^3};\Z )\overset {i ^*}\to H^4(M;\Z )\end{equation}
is exact. Next consider the bundle $E_{S^3}\to M_{S^3}$ associated to 
the $S^3$-equiva\-riant vector bundle $E\to M$ via the Borel 
construction. We denote its first Pontrjagin class by $p_1(E)_{S^3}\in 
H^4(M_{S^3}; \Z )$. Since $i^*(p_1(E)_{S^3})=p_1(E)=0$ it follows from 
(\ref{eqspecseq}) that $p_1(E)_{S^3}=\pi^*(\alpha)$ for some $\alpha 
\in H^4(BS^3;\Z )$. By naturality, 
$p_1(E)_{S^1}=\pi ^*(C^{\prime \prime }\cdot x^2)$, where $\pi $ also 
denotes $M_{S^1}\to BS^1$, $x\in H^2(BS^1;\Z )$ denotes a generator 
and $C^{\prime 
\prime }\cdot x^2$ is the image of $\alpha $ under the map induced by 
the inclusion $S^1\hookrightarrow S^3$. Note that $p_1(E)_{S^1}$ 
reduces at $q\in Y_s$ to $(\sum 
_{i=1}^r e_{s,i}^2)\cdot x^2\in H^4(BS^1;\Z )\cong H^4(q_{S^1};\Z )$. Hence $\sum 
_{i=1}^r e_{s,i}^2=C^{\prime \prime }$ for any $Y_s$.
\proofend

\section{Bordism-finiteness}\label{bordfin}
In this section we apply the \LFF \ to the result of the previous 
section to derive bordism-finiteness theorems. Let $\partial $ be an 
$S^1$-equivariant elliptic operator over $M$ with symbol $\sigma \in 
K_{S^1}(TM)$. By the fundamental work of Atiyah and Singer (cf. 
\cite{AtSiI}) the equivariant index of $\partial $ is equal to the 
topological index $\tind (\sigma )$, where $\tind :K_{S^1}(TM)\to 
R(S^1)$ is defined as the push-forward in complex $K$-theory for $M\to 
pt$. Let $N$ denote the normal bundle of $i:M^{S^1}\hookrightarrow M$ 
and $\Lambda 
_t:=\sum_{i=0}^\infty \Lambda ^i\cdot t^i$. 
\begin{theorem}[\LFF \ \cite{AtSe}]\label{lefschetz}
For an element $\sigma \in K_{S^1}(TM)$ and any topological generator 
$\lambda 
\in S^1$ \begin{equation}\label{eqLFF} 
\tind (\sigma 
)(\lambda )=\tind 
\left(\frac {i^*(\sigma )(\lambda )}{\Lambda _{-1}(N\otimes _\R \C 
)(\lambda)}\right ). 
\end{equation}
\proofend
\end{theorem}

\bigskip
\noindent 
Note that the right hand side of formula (\ref{eqLFF}) consists of a 
finite sum of local contributions at the fixed point components of the 
$S^1$-action. Since the set of topological generators is dense in 
$S^1$ the theorem above implies that the index function 
$\tind:K_{S^1}(TM)\to R(S^1)$ vanishes identically if $M^{S^1}$ is 
empty. 

If the elliptic operator is geometrically defined using an 
$H$-structure (cf. \cite{AtSiIII}), for example if $\partial$ is a 
twisted signature operator, then the local data in formula 
(\ref{eqLFF}) may be expressed in terms of the tangent bundle and the 
vector bundles associated to the $H$-structure restricted to 
$M^{S^1}$.

Next recall that the Pontrjagin numbers of $M$ are determined by 
twisted signatures, where the twist bundles are associated to the 
tangent bundle via some orthogonal representation. Since the 
$S^1$-action lifts canonically to these bundles one can consider the 
$S^1$-equivariant twisted signatures. If the $S^1$-action on $M$ has 
no fixed point these indices have to vanish since the index function 
$\tind$ vanishes identically. In this case all Pontrjagin numbers of 
$M$ vanish (for an elementary proof see \cite{Bo}) which implies that 
$M$ represents an element of order two in the oriented bordism ring. 

The following reformulation is convenient for our purposes. Given two 
$S^1$-manifolds $M_1$ and $M_2$ we say that they have the same {\bf local $\mathbf{S^1}$-geometry} if there exists an $S^1$-equivariant 
orientation preserving diffeomorphism from the normal bundle of the 
fixed point manifold $M^{S^1}_1$ to the normal bundle of $M^{S^1}_2$. 
In this case one can glue $M_1$ and $-M_2$ along the fixed point 
manifold together to obtain a manifold with fixed point free 
$S^1$-action which is bordant to $M_1 
-M_2$ (cf. for example \cite{CoFl}, Th. (22.1)). We summarize the 
discussion in 
\begin{corollary}\label{localgeometry} Let $M_1$ and $M_2$ be
$S^1$-\mfd s with the same local $S^1$-geometry. Then $M_1 -M_2$ 
represents a torsion element in the oriented bordism ring.\proofend 
\end{corollary}

\bigskip
\noindent
We shall combine the corollary with Lemma \ref{mainlemma} to obtain 
bordism-finiteness for the class of manifolds for which a negative 
multiple of the first Pontrjagin class is a sum of squares. Let us 
call such manifolds {\bf semi-negative}, i.e. a manifold $M$ is 
semi-negative if there exists a finite number of classes $y_j\in 
H^2(M;\Z )$ and a positive integer $N$ such that $N\cdot p_1(M)+\sum 
y_j ^2=0$. Note that $M$ is semi-negative if and only if $M$ admits 
complex line bundles $L_j$ for which the first Pontrjagin class of 
$E:=\sum L_j+N\cdot TM$ vanishes. Of course many manifolds such as 
$3$-connected manifolds with non-vanishing first rational Pontrjagin 
class (e.g. quaternionic projective space $\H P^k$, $k>1$) cannot be 
semi-negative. However the class of semi-negative manifolds is quite 
rich and contains some interesting families. 
\begin{remarks}\label{s-nex} The class of semi-negative manifolds contains
\begin{enumerate}
\item manifolds with $p_1$ torsion, e.g. $BO\langle 8\rangle $-manifolds or $4$-connected manifolds.
\item manifolds with cohomology ring in degree $\leq 4$ like a $4$-manifold with indefinite intersection form.
\item in a given dimension all but a finite number of complete intersections.
\item infinitely many homotopy complex projective spaces in any given dimension $2n\geq 6$.
\end{enumerate}
\end{remarks}

\begin{theorem}\label{theobordfin} Let $C$ and $m$ be natural numbers. The class 
of semi-negative connected $m$-dimensional oriented manifolds with Euler characteristic 
$\leq C$ and an $S^3$-action satisfying $(\ast )$ contains only 
finitely many oriented bordism types. 
\end{theorem}

\bigskip
\noindent
{\bf Proof:} The theorem follows from Corollary \ref{localgeometry} 
once we know that the local $S^1$-geometry is determined by $(m,C)$ up 
to finite ambiguity. Since the oriented bordism group $\Omega 
^{SO}_m$ is finite in dimension $m\not 
\equiv 0 \bmod 4$ we may assume $m=4l$.  Let $M$ be a $4l$-dimensional manifold with 
$S^3$-action satisfying $(\ast )$ and Euler characteristic $\chi 
(M)=C^\prime 
\leq C$. By the classical \LFF \ the induced $S^1$-action has 
$C^\prime $ isolated fixed points. Next assume that $M$ is 
semi-negative, i.e. $p_1(L_1+\ldots +L_k+N\cdot TM)=0$ for a positive 
integer $N$ and complex line bundles $L_j$. By Lemma \ref{mainlemma} 
the tangential weights of $M$ are determined by the dimension of $M$ 
up to finite ambiguity. So all manifolds in the theorem with Euler 
characteristic equal to $C^\prime $ represent only a finite set of 
local $S^1$-geometries. Since $0\leq C^\prime \leq C$ the theorem 
follows from Corollary \ref{localgeometry}.\proofend 

\bigskip
\noindent
The theorem immediately implies Theorem \ref{p1zero}. It extends to 
the complex bordism ring if one assumes in addition that the induced 
$S^1$-action preserves the stable almost complex structure. A 
corresponding result for the $Spin^c$-bordism ring does not hold since 
an $S^3$-action always lifts to a given $Spin^c$-structure and there 
are just too many of them. However Theorem \ref{theobordfin} admits 
the following refinement. 

Let $X(k)$ be the Cartesian product of $BSO$ and $k$ copies of $\C 
P^\infty $. For a fixed positive integer $N$ let $f_N:X(k) 
\to K(\Z ,4)$ be the map which classifies $-N\cdot p_1+\sum_{j=1}^kx_j^2\in H^4(X(k);\Z )$. Here $p_1\in H^4(BSO ;\Z )$ denotes 
the universal first Pontrjagin class and $x_j$ denotes a generator for 
the $j$-th copy of $H^2(\C P ^\infty ;\Z )$. Let $B(k,N)\to X(k)$ be 
the pullback of the path fibration $E(\Z ,4)\to K(\Z ,4)$ via $f_N$ 
and let $\pi 
:B(k,N)\to BSO $ denote the projection. Fix a classifying map $\nu :M\to BSO$ for the stable 
normal bundle of $M$. 

A $B(k,N)$-structure for $M$ is the isotopy class $h$ of a lift $M\to 
B(k,N)$ of $\nu$ in the fibration $\pi :B(k,N)\to BSO$. For fixed $k$ 
and $N$ we call $(M,h)$ a {\bf $\mathbf B$-manifold}. Note that 
$(M,h)$ determines $k$ complex line bundles $L_j$ over $M$ such that 
$N\cdot p_1(M)+\sum_{j=1}^k p_1(L_j)=0$. In particular, $M$ is 
semi-negative. 

Any polynomial in $p_i(M)$ and $c_1(L_j)$ is called a {\bf 
characteristic class} of the $B$-manifold $(M,h)$. The {\bf 
characteristic numbers} are defined as the values of the 
characteristic classes on the fundamental cycle of $M$. By the 
Pontrjagin lemma the characteristic numbers only depend on the bordism 
class of $(M,h)$. The bordism group $\Omega _*^{B}$ of $B$-manifolds 
may be studied in terms of stable homotopy theory using the 
Pontrjagin-Thom construction. 
\begin{proposition}\label{bkstructure} For fixed $k$ and $N$ the bordism group $\Omega 
_n^{B}$, $n\in \Z $,
 is finitely generated. The group 
$\Omega _n^{B}\otimes \Q $ is completely determined by the 
characteristic numbers. 
\end{proposition}

\bigskip
\noindent
{\bf Proof:} Since the argument is standard (cf. \cite{La}, \cite{St}) 
we only sketch it. First apply the construction above to $BSO(r)$ for 
$r\in \N $. So $X_r$ is the product of $BSO(r)$ and $k$ copies of $\C 
P 
^\infty $, $B_r\to X_r$ is obtained as pullback of the path fibration $E(\Z ,4)\to K(\Z ,4)$ and $\pi 
_r:B_r\to BSO(r)$ is the projection. Note that $B(k,N)
=\lim _{r\to \infty} B_r$. Let $\gamma _r\to B_r$
 be the pullback of the 
universal vector bundle over $BSO(r)$ via $\pi _r$ and let $M(\gamma 
_r)$ be its Thom space. The Pontrjagin-Thom construction gives an isomorphism $\Phi :\Omega 
_n^{B}\to \pi _{n+r}(M(\gamma 
_r))$ for $r\gg n$, which defines an isomorphism
$$\Phi :\Omega _n^{B}\to  \lim _{r\to \infty} \pi _{n+r}(M(\gamma 
_r))\cong \pi _n(M(\gamma ))$$ between the bordism group of $n$-dimensional $B$-manifolds 
and the $n$-th homotopy of the associated spectrum $M(\gamma )$. Since 
$M(\gamma 
_r)$ is a CW-complex with finite skeletons $\Omega _n^{B}$ is finitely generated.

Next assume the $n$-dimensional $B$-manifold $(M,h)$ has vanishing 
characteristic numbers. For the second statement it suffices to show 
that the bordism class of $(M,h)$ vanishes in $\Omega ^B_n\otimes 
\Q $. Note that 
the space of characteristic classes of $(M,h)$ is equal to $\hat 
h^*(H^*(B_r;\Q ))$, $r\gg n$, where $\hat h:M\to B_r$ represents $h$. 
Since $(M,h)$ has vanishing characteristic numbers $\hat h_*(\mu _M)$ 
vanishes in $H_n(B_r;\Q )$. Now apply the Thom isomorphism for the 
normal bundle of $M$ and the bundle $\gamma _r\to B_r$ to conclude 
that the composition of the Pontrjagin-Thom map $\Phi $ and the 
rational Hurewicz homomorphism $\Psi $ 
$$\Omega 
_n^{B}\otimes 
\Q \overset {\Phi }\to \pi _{n+r}(M(\gamma _r))\otimes \Q  
\overset {\Psi }\to H_{n+r}(M(\gamma _r);\Q ),$$ $r\gg n$, maps $(M,h)$ to zero. Since $\Phi $ and
 $\Psi $ are isomorphisms $(M,h)$ vanishes in $\Omega ^B_n\otimes 
\Q $.\proofend 

\begin{theorem}\label{bkfinite} Let $C$ and $m$ be natural numbers. For fixed $k$ and $N$ the class 
of connected $m$-dimensional $B$-manifolds with Euler characteristic $\leq C$ 
and an $S^3$-action satisfying $(\ast )$ contains only finitely many 
$B$-bordism types. 
\end{theorem}

\bigskip
\noindent
{\bf Proof:} By Proposition \ref{bkstructure} we may assume that 
$m=2n$. Let $(M,h)$ be a $2n$-dimensional $B$-manifold with Euler 
characteristic $\leq C$ and an $S^3$-action satisfying $(\ast )$. The 
map $h$ induces via the projection of $B$ to the $k$-fold product of 
$\C P^\infty $ a classifying map for $k$ complex line bundles $L_j$ 
which satisfy $N\cdot p_1(M)+\sum_{j=1}^k p_1(L_j)=0$. 

We want to show that the characteristic numbers of $(M,h)$ are 
determined by $(m,C,k,N)$ up to finite ambiguity. To this end let $J$ 
denote the subring of $K(M)$ generated by the complex line bundles 
$L_j$ and vector bundles associated to the tangent bundle. Let 
$ch:K(M)\to H^*(M;\Q )$ be the Chern character. We note that 
$ch(J)\otimes 
\Q $ is the subspace $V$ of $H^*(M;\Q )$ which is 
spanned by the characteristic classes of the $B$-manifold $(M,h)$. 

Next we identify the characteristic numbers of $(M,h)$ with certain 
twisted signatures. By the cohomological version of the index theorem 
(cf. \cite{AtSiIII}) the index of the signature operator of $M$ 
twisted with $F\in J$ is given by 
$$sign(M;F)=\left\langle \prod _{i=1}^n(u_i\cdot \frac {1+e^{-u_i}}{1-e^{-u_i}})
\cdot ch (F),\mu _M\right \rangle ,$$
where $\pm u_i$ are the formal roots of $M$, $\mu 
_M$ is the fundamental cycle of $M$ and $\langle \quad ,\quad \rangle $
 denotes the pairing between cohomology and homology. This implies that 
$\{sign(M;F)\, \mid \, F\in J\}$ spans the $\Q $-vector space $\langle 
V,\mu _M\rangle $. Hence the characteristic numbers of $(M,h)$ are 
determined by twisted signatures, where the twist bundle is an element 
of $J$. 

We are now in the position to prove the theorem from the \LFF . Choose 
a finite set $\{ E_1,\ldots ,E_r\}\subset J$ ($r$ depends on $k$ and 
the dimension of $M$) of vector bundles which span $J\otimes \Q $. 
Here each $E_i$ is given by a universal polynomial (which only depends 
on $(m,k)$) in the complex line bundles $L_j$ and vector bundles 
associated to the tangent bundle. We view $E_i$ as an 
$S^3$-equivariant vector bundle by lifting the $S^3$-action (uniquely) 
to each $L_j$. By Lemma \ref{mainlemma} the tangential weights and the 
weights of $L_j$ are determined by $(m,N)$ up to finite ambiguity. 
This implies the same for the weights of $E_i$. 

Next consider the signature operator twisted by $E_i$. By the 
discussion of the \LFF \ after Theorem \ref{lefschetz} its 
$S^1$-equivariant index is equal to a sum of $\leq C$ local 
contributions which only depend on the weights of $E_i$ and the 
tangential weights. Since these belong to a finite set which only 
depends on $(m,k,N)$ we conclude that the ordinary twisted signatures 
$sign(M;E_i)$, $i=1,\ldots , r$, are determined by $(m,C,k,N)$ up to 
finite ambiguity. This implies the same for the characteristic numbers 
of the $B$-manifold $(M,h)$. Now the theorem follows from Proposition 
\ref{bkstructure}.\proofend

\section{Examples of semi-simple group actions}\label{examples}

In this section we show that Theorem \ref{theobordfin} is sharp in the 
sense that bordism-finiteness fails if one weakens the assumptions on 
the $S^3$-action or the first Pontrjagin class. As explained in the 
introduction one cannot expect bordism-finiteness for $S^1$-actions if 
one allows arbitrary fixed point sets. We restrict to S$^1$-actions 
with isolated fixed points which satisfy a prescribed upper bound for 
the number of fixed points. Also we assume that the action extends to 
an action of $S^3$. In this situation, by Theorem \ref{theobordfin}, 
bordism-finiteness holds if the $S^3$-action has a fixed point and the 
manifolds are semi-negative. The following two propositions 
show that both assumptions are necessary. 

\begin{proposition}\label{firstex} There exist connected $20$-dimensional semi-nega\-tive manifolds $M_l$, $l\in \N $, with Euler 
characteristic equal to $12$ 
 which represent distinct oriented bordism classes such that each $M_l$ supports a fixed 
point free $S^3$-action with isolated $S^1$-fixed points. 
\end{proposition}
\begin{proposition}\label{secondex} There exist connected $20$-dimensional manifolds $N_l$, $l\in \N $, with Euler 
characteristic equal to $23$ which represent distinct oriented bordism 
classes such that each $N_l$ supports an $S^3$-action with fixed point 
and 
  isolated $S^1$-fixed points.
\end{proposition}

\bigskip
\noindent
Note that by Theorem \ref{theobordfin} almost all of the $N_l$ are not 
semi-negative. We remark that examples as in the propositions are 
necessarily of dimension $4k\geq 8$ since in dimension $4$ the upper 
bound on the number of isolated fixed points gives a bound on the absolute value of the 
signature. The remaining part of this section is devoted to the proof 
of the propositions above. We use a kind of induction to extend 
actions on the base of a fibre bundle to the total space if the 
fibration is associated to a principal torus bundle. 

Let $G$ be a compact connected Lie group which acts from the left on a 
connected manifold $Z$. Let $F$ be a connected manifold with left 
$U(1)$-action. Assume the first Betti number $b_1(Z)$ vanishes or $G$ 
is simply-connected. 

For $y\in H^2(Z;\Z )$ let $S\to Z$ denote the $U(1)$-\pb \ with first 
Chern class equal to $y$ and let $M:=S\times _{U(1)}F$ be the 
associated fibre bundle. In \cite{HaYo} it was shown that the 
$G$-action lifts to $S$ (in fact uniquely if $G$ is simply-connected). 
For a fixed lift $G$ acts on $M$ by $g(s,f)_\sim :=(gs,f)_\sim $. Note that the projection $M\to Z$ is $G$-equivariant. 

We now restrict to the case $G=S^3$ and fix $S^1\hookrightarrow S^3$. 
Let $a_Y$ denote the weight of the $S^1$-action on $S$ restricted to a 
connected component $Y$ of $Z^{S^1}$, i.e. at a point of $Y$ (and 
hence at any point of $Y$) the $S^1$-action on the fibre of $S$ has 
character $\lambda \mapsto 
\lambda ^{a_Y}$. For further reference we note the elementary
\begin{lemma}\label{fixed points} Assume none of the weights $a_Y$ vanish. Then the $S^1$-action on $M$ has isolated fixed 
points if this holds for the $S^1$-action on $Z$ and the $U(1)$-action 
on $F$.\proofend 
\end{lemma}

\noindent 
We now begin with the construction of the examples mentioned before. 
The first series consists of Cayley plane bundles over $Z:=S^2\times 
S^2$ which support an $S^3$-action with isolated $S^1$-fixed points 
but no $S^3$-fixed point. We take the left $S^3$-action on $Z$ induced 
from the homogeneous action on each copy of $S^2$. This action has 
isolated $S^1$-fixed points but no $S^3$-fixed points. Its principal 
stabilizer is equal to $\Z _2$, the center of $S^3$. 

As fibre $F$ we take the Cayley plane $Cl_2=F_4/Spin(9)$. We fix an 
orientation of $Cl_2$. Next we choose an embedding 
 $j:U(1)\hookrightarrow T$, where $T$ is a maximal 
torus of $Spin(9)$, such that the induced $U(1)$-action on $Cl_2$ is 
effective and has isolated fixed points. 

Let $\gamma _i\to Z$ denote the pullback of the Hopf bundle over $S^2$ 
under the projection of $Z$ on the $i$-th factor and let 
$z_i:=c_1(\gamma _i)\in H^2(Z;\Z )$. Equip $Z$ with the orientation 
dual to $z_1\cdot z_2$. Note that the weights of the $S^3$-equivariant 
line bundle $\gamma 
_1^{a}\otimes \gamma _2^{b}$ at the $S^1$-fixed points of $Z$ 
have the form $\pm a\pm b$. 

To construct $M_l$ fix positive integers $a\not \equiv b\bmod 2$, 
let $S_l$, $l\in \N $, be the $U(1)$-principal bundle associated to 
$(\gamma _1^a\otimes \gamma 
_2^{b})^{2l+1}$ and set $M_l:=S_l\times 
_{U(1)}Cl_2$. Note that $M_l$ comes with 
an orientation induced from the orientations of $Z$ and $Cl_2$. As 
explained above the $S^3$-action on $Z$ lifts to $M_l$. For the 
induced $S^1$-action some of the tangential weights of $M_l$ are odd 
since the weights of $S_l$ have the form $(2l+1)\cdot (\pm a\pm b)$ 
and the action of $U(1)$ on $Cl_2$ is effective. In particular, the 
principal isotropy group of the $S^3$-manifold $M_l$ is trivial. 

\bigskip
\noindent
{\bf Proof of Proposition \ref{firstex}:} The Euler characteristic 
$\chi (M_l)$ of the oriented $20$-dimensional manifold $M_l$ is just the 
product of the Euler characteristic of the base and the fibre, thus 
$\chi (M_l)=12$. The $S^3$-action on $M_l$ has no fixed points since 
$p:M_l\to Z$ is equivariant. By Lemma \ref{fixed points} the induced 
$S^1$-action on $M_l$ has isolated fixed points. Since $H^*(M_l;\Z )$ 
is isomorphic to $H^*(Z;\Z )$ in degree $\leq 4$ (the fibre is 
$7$-connected) each $M_l$ is semi-negative (see Remarks \ref{s-nex}). 

We compute the Milnor number $\langle s_{10}(TM_l),\mu 
_{M_l}\rangle $ to show that the manifolds $M_l$ represent distinct bordism 
classes (for a real vector bundle with formal roots $\pm u_i$ the 
class $s_{2t}$ is defined as $\sum u_i^{2t}$). Let $\pi: E\to BU(1)$ 
denote the pullback of $BSpin(9)\to BF_4$ to $BU(1)$ via $Bj$ and let 
$E^\triangle 
\to E$ denote the tangent bundle along the fibres. One computes using 
\cite{BoHi} that $\pi 
_!(s_{10}(E^\triangle ))=\alpha \cdot 
x^2$, where $\alpha \neq 0$ (cf. for example \cite{DDISS}, Section 
4.2). Here $x$ is a generator of $H^2(BU(1);\Z )$ and $\pi _!$ is the 
push-forward in cohomology. Next note that the map $f:Z\to BU(1)$ 
which classifies $(2l+1)\cdot (a\cdot z_1+b\cdot z_2)$ is covered by 
$\tilde f:M_l\to E$ and the tangent bundle along the fibres of 
$p:M_l\to Z$ is isomorphic to $\tilde f^*(E^\triangle )$. Now compute 
$$\langle s_{10}(TM_l),\mu _{M_l}\rangle =
\langle s_{10}(\tilde f^*(E^\triangle )\oplus p^*(TZ)),\mu _{M_l}\rangle =
\langle p_!(s_{10}(\tilde f^*(E^\triangle ))),\mu _Z\rangle =$$
$$\langle f^*(\pi _!(s_{10}(E^\triangle ))),\mu _Z\rangle =
\langle \alpha \cdot (2l+1)^2\cdot (a\cdot z_1+b\cdot z_2)^2,\mu_Z\rangle =
2\cdot \alpha \cdot (2l+1)^2\cdot a\cdot b.$$ The computation shows 
that the $M_l$ represent distinct bordism classes. 
\proofend

\bigskip
\noindent The next series is constructed from the series above. Equip $\C P^{10}$ with the
$S^3$-action induced by the direct sum of the trivial one-dimensional 
complex representation and the irreducible complex 
$S^3$-representation of dimension $10$. Note that $S^3$ acts on $\C 
P^{10}$ with fixed point, the principal isotropy group is trivial and 
the induced $S^1$-action has isolated fixed points. The principal 
isotropy group of $M_l$ is also trivial (see the discussion before the 
proof of Prop. \ref{firstex}). Now define $N_l$ by taking the 
equivariant connected sum of $M_l$ and $\C P^{10}$ along a principal 
orbit. 

\bigskip
\noindent {\bf Proof of Proposition \ref{secondex}:} First note that the disjoint
union of $M_l$ and $\C P^{10}$ is bordant to $N_l$. Also $\chi (N_l)=\chi (M_l)+\chi (\C P^{10})$. Hence, by Proposition \ref{firstex} the manifolds $N_l$, $l\in 
\N $, represent distinct bordism classes and have Euler characteristic $\chi (N_l)=23$.
 By construction $S^3$ acts on $N_l$ 
with fixed point and the induced $S^1$-action has isolated fixed 
points.\proofend 

\section{Homotopy complex projective spaces and complete intersections}\label{homcp}

In this section we prove the theorem on homotopy complex projective 
spaces given in the introduction and a related result for complete 
intersections (see Theorem \ref{complete} below).

Let $M$ be a $2n$-dimensional closed manifold with $H^*(M;\Z )\cong 
H^*(\C P^n;\Z )$ and non-trivial $S^1$-action. Let $\gamma $ denote 
the complex line bundle with $c_1(\gamma )=x$, where $x$ is a fixed 
generator of $H^2(M;\Z )$. By \cite{HaYo} the $S^1$-action lifts to 
$\gamma $. We denote the weight of $\gamma $ and the tangential 
weights at a connected component $Y_s$ of $M^{S^1}$ by $a_s$ and 
$m_{s,1},\ldots ,m_{s,n}$, respectively. In the proof we use the 
following information on the weights: For fixed $s$ the weights of $\gamma $ are related to the tangential weights 
by 
\begin{equation}
\label{cpweight} \vert \prod_{t\neq s}(a_t-a_s)^{n_t+1}\vert =\vert 
\prod_{m_{s,i}\neq 0} m_{s,i}\vert,
\end{equation}
where $n_t$ denotes the complex dimension of $Y_t$. To prove this 
identity one either applies the localization theorem in $K$-theory to 
the $S^1$-action and induced $\Z 
_p$-actions for $p$ a large prime number (cf. \cite{Pe}, Th. 2.8) or 
uses cohomological means (cf. \cite{Br}, Ch. VII, Th. 5.5). 

\bigskip
\noindent
{\bf Proof of Theorem \ref{hcp}:} For $n=1$ the theorem is trivial. For $n=2$ it follows for example from the classification of $4$-dimensional $S^3$-manifolds (cf. \cite{MePa}). So assume $n\geq 3$.

Let $M$ be an oriented homotopy $\C P^n$ with $S^3$-action satisfying 
$(\ast)$. Let $Y_0,\ldots ,Y_n$ denote the isolated fixed points under 
the induced $S^1$-action on $M$. We assume that $Y_0$ is also fixed by 
$S^3$. Hence, $a_0$ vanishes and the tangential weights $m_{0,1}, 
\ldots ,m_{0,n}$ are determined by the dimension $m=2n$ up to finite 
ambiguity. Therefore $\vert \prod_{i=1}^n m_{0,i}\vert $ is 
 bounded from above by a constant which only 
depends on $m$. Now apply formula (\ref{cpweight}) twice to conclude that the same holds for the 
absolute value of all the weights $a_s$ and $m_{s,i}$. 

Next we argue as in the proof of Theorem \ref{bkfinite} to show
that the values of polynomials in the Pontrjagin classes and the 
generator $x\in H^2(M;\Z )$ on the fundamental cycle belong to a 
finite set which only depends on $m$. Since the cohomology ring of $M$ 
is a truncated polynomial ring in $x$ the Pontrjagin classes are also 
determined up to finite ambiguity. By simply-connected surgery theory 
it follows that the diffeomorphism type of $M$ belongs to a 
finite set which only depends on the dimension $m$. 
\proofend

\bigskip
\noindent
Next we consider complete intersections. A complete intersection 
$V_n^{(d_1,\ldots ,d_r)}$ of complex dimension $n$ and multidegree 
$(d_1,\ldots ,d_r)$, $d_i\geq 2$, is defined as the transversal 
intersection of hypersurfaces of degree $d_i$, $i=1,\ldots ,r$, in $\C 
P^{n+r}$. Thom showed that the diffeomorphism type of $V_n^{(d_1,\ldots 
,d_r)}$ is completely determined by $n$ and $(d_1,\ldots ,d_r)$. 
\begin{theorem}\label{complete} In a fixed complex dimension there are only finitely many
complete intersections with an $S^3$-action satisfying $(\ast 
)$.\end{theorem} 

\noindent
We remark that the $\hat A$-vanishing theorem of Atiyah-Hirzebruch 
(cf. \cite{AtHi}) may be used to show that there are only finitely 
many $Spin$-complete intersections with non-trivial $S^1$-action in a fixed even complex 
dimension. Also results of Hattori (cf. \cite{Ha}) on the equivariant 
$Spin^c$-Dirac operator imply finiteness of the number of complete 
intersections with non-trivial $S^1$-action if one assumes that the 
action preserves the induced stable almost complex structure. 

\bigskip
\noindent {\bf Proof of Theorem \ref{complete}:} Assume $M:=V_n^{(d_1,\ldots ,d_r)}$ admits an
$S^3$-action satisfying $(\ast )$. Let $\gamma $ denote the pullback 
of the dual Hopf bundle over $\C P^{n+r}$ via $i:M\hookrightarrow \C 
P^{n+r}$ and let $x:=c_1(\gamma )\in H^2(M;\Z )$. Recall from 
\cite{Hi} that $p(M)=(1+x^2)^{n+r+1}\cdot 
\prod 
_{i=1}^r(1+d_i^2\cdot x^2)^{-1}$ and $\langle x^n,\mu _M\rangle =\prod _{i=1}^r d_i$. Hence, for all but a finite number of 
multidegrees $M$ is semi-negative with $N=1$ and $L_j=\gamma $. For 
semi-negative $M$ the tangential weights and the weights of $\gamma $ 
are determined by the complex dimension $n$ up to finite ambiguity by 
Lemma \ref{mainlemma}. We argue as before (see the proof of Th. \ref{bkfinite}) to conclude that the values of polynomials in the 
Pontrjagin classes and $x$ on the 
fundamental cycle $\mu_M$ belong to a finite set which only depends on 
$n$. In particular this holds for $\langle x^n,\mu_M\rangle =\prod 
d_i$. Since $d_i\geq 2$ the theorem follows.\proofend 

\baselineskip=3mm
\setlength{\parindent}{0pt}

{\small

\bigskip
\noindent
Anand Dessai\\
Department of Mathematics\\
University of Augsburg\\
D-86159 Augsburg\\
e-mail: dessai@math.uni-augsburg.de\\
http://www.math.uni-augsburg.de/geo/dessai.html

\end{document}